\documentclass{amsart}
\usepackage{times}
\usepackage{amssymb,latexsym}
\usepackage[T1]{fontenc}
\usepackage[english]{babel}
\usepackage{graphicx}

\theoremstyle{definition}

\theoremstyle{definition}

\theoremstyle{definition}

\theoremstyle{definition}

\title{To found or not to found? That is the question!}
\author{davide bondoni}
\address{via Bersaglio 2\\
25070 - Anfo (BS)\\
Italy}
\email{davidebond@yahoo.it}
\urladdr{http://www.davidebondoni.eu}

\keywords{Dedekind, Schr\"oder, principle of induction, calculus of relatives, theory of chains}
\subjclass[2010]{Primary: 01A55; Secondary: 03D70}

\begin{document}

\begin{abstract}
Aim of this paper is to confute two views, the first about Schr\"oder's presumptive \emph{foundationalism}, according to he founded mathematics on the calculus of relatives; the second one mantaining that Schr\"oder only in his last years (from 1890 onwards) focused on an universal and symbolic language (by him called {\it pasigraphy}). We will argue that, on the one hand Schr\"oder considered the problem of founding mathematics already solved by Dedekind, limiting himself in a mere translation of the Chain Theory in the language of the relatives. On the other hand, we will show that Schr\"oder's pasigraphy was connaturate to himself and that it roots in his very childhood and in his love for foreign languages.
\end{abstract}

\maketitle

\section{Introduction}

The present article develops in two sections: the first one (Section \ref{found}), devoted to refute the opinion that Schr\"oder founded mathematics on the calculus of relations, the second one (Section \ref{found2}), devoted to prove that Schr\"oder already in his youth was interested in the various forms of language. As the the reader can see, there is an overlapping between these parts, as the calculus of relations was choosen by Schr\"oder as a suitable universal language in which express mathematics. In this sense, the calculus of relations must be regarded on the background of other similar efforts, as Peano's \emph{lingua franca}, aiming to find a purely symbolic language for mathematics. I know that Schr\"oder states explicitely that the main goal of his calculus of relatives is to give a definition of number, but, as argued elsewhere, I believe that it would be more appropriate saying that Schr\"oder tried to translate in terms of relations Dedekind's definition of number.\footnote{Schr\"oder seems not be aware that Dedekind did not give a definition of \emph{number}, but the definition of \emph{set of numbers}.}\\
\indent Some scholars may rebuke my interpretations, invoking the faithfulness to Schr\"oder's own words. To them, I reply that sometimes we need to go beyond the literal expressions to grasp the meaning in question, more or less as a psychiatrist does in analyzing the not-said, the unconscious.\footnote{I usually make the following example: in \cite[p.~322]{schr69} Schr\"oder writes: \emph{(\ldots) f\"ur alle Punkte $z$ eines gewissen um den Punkt $z_1$ herum liegende Gebiet strebt die ohne Ende fort iterirte Funktion $F(z)$ der Wurzel $z_1$ der Gleichung $F(z) = z$ als Grenze zu}. Alexander \cite[p.6]{alex} translates \emph{Gebiet} with \emph{area}. This is linguistic correct, but it hampers a true comprehension, beacause Schr\"oder means with 'Gebiet' \emph{Umgebung} [i.e. neighbourhood] as it is evident from the context.} But it is not only a matter of interpretation. Schr\"oder is ambigous on this point: if in the first lectures of the third volume of the \emph{Vorlesungen \"uber die Algebra der Logik} he seems interested in founding mathematics on the calculus of relations, with the passing of the time, he become aware of the power of the theory of relation as a symbolic language. From that point onwards, Schr\"oder puts aside the calculus of relatives in itself to translate set-theoretic problems in his new language. Now Schr\"oder recognizes that his calculus of relations is not only a calculus but it is also a language: a language more expressive that Peano's one, for example.\\
\indent This is my rationale to support the thesis that from the third volume of the \emph{Vorlesungen}, devoted to the calculus of relations, Schr\"oder turned to set-theoretic problem, as that concerning the well-foundness \cite{schr01}. On this point too there is no agreement between the scholars. They insist that Schr\"oder was interested in relations in themselves, but that this work was interrupted by a repentine death. Nevertheless, from the publication of the third volume of the \emph{Vorlseungen} (1895) and Schr\"oder's death (1902) passed seven year that he could spend in many ways. Neither it can be said that Schr\"oder was diverted by his university duties, because he was rector of Karlsrue university one year only, or that he was unable to copy with his numerous hobbies, as Dipert mantains:
\begin{quotation}
Among Schr\"oder's hobbies were hiking, swimming, ice-skating, horseback riding, and gardening\footnote{In 2010 I discovered a little note by Schr\"oder on grafting \cite{schr88}.} -- and perhaps these are what L\"uroth is also suggesting distracted Schr\"oder from his research.\footnote{\cite[p.~126]{dipert2}.}
\end{quotation}

Schr\"oder had the possibility to focus on one matter in place of another. If the German mathematician would continue in investigating the calculus of relations in itself, and not as a language in which tackling set-theoretical questions, he had the free will to do it. As a matter of fact, he was more allured by set-theory. It is not so difficult to understand. As a human being, Schr\"oder had the choice to think what he more liked. \emph{Die Gedanken sind frei} sings the Mahler prisoner in the tower. Albeit condemned to death, he is free to think what he will.\\
\indent Summarizing: the major part of scholars engaged in Schr\"oder consider him as a logician; i.e. a mathematician that from 1877 converted his activity to logic. I confute this view. Schr\"oder was and remained all along his life a mathematician. He searched for new fields of work, but not leaving a mathematical style of thought. For this reason, he passed from the calculus of relations to the more appealing (and mathematical) set-theory.\footnote{To my present knowledge, no one historian of logic pointed out that in the fifth lecture of the \emph{Vorlesungen} Schr\"oder searches for fixed point results.}\\
\indent Finally, Schr\"oder's pasigraphy and Schr\"oder`s theory of relations are two faces of the same object, that we can see from different perspectives. It is right asserting that from a point of view and in a precise lapse of time, Schr\"oder regarded the theory of relations from a computational point of view, but it is also right asserting that in another time he regarded this theory as a symbolic language. This must be conceded. Then, my article is devoted to the same theory viewed in two different ways in two different moments.\\
\indent With the following section I ponder on Schr\"oder's presumptive foundationalism, casting light on the time when Schr\"oder passed to consider the theory of relations a language in which \emph{translating} and not \emph{founding} mathematical concepts.

\section{To found: that is the question}\label{found}
Elsewhere I asserted that Ernst Schr\"oder was to be considered a mathematician and not a logician, being his knowledge in issue very meager, and because he was engaged all life long in mathematical questions (as in the algebraic {\it Solution Problem}).\footnote{See \cite{algebra} and \cite{gruppi}. Of the same opinion was Randall R.~Dipert, asserting: {\it 
Schr\"oder was a \textbf{practicing} mathematician after all, and his influence on philosophical discussion, other than indirectly through later mathematical logic, seems to have been very small} \cite[p.~140, the emboldening is mine]{dipert}.} Schr\"oder had a {\it structural} view of mathematics, according to, any concept has a meaning only in virtue of the place it takes inside the theory.\footnote{More on this topic in \cite{LL1} and \cite{LL2}.} In other words, the mathematical concepts are context-dependent, being this context a relational lattice.\footnote{I will use the word 'lattice' in a non technical way.} Well, if a mathematical theory is only a set of relations, what does Schr\"oder mean with {\it relation}?\\
\indent This question lays at the core of Schr\"oder's work, because mathematical formulas are only strings of symbols without interpretation. For Schr\"oder, a formula means {\it something} only inside the theory in which is formulated, and such theory is a structured set of relations.\footnote{This explain Schr\"oder interest in what we call today {\it Group Theory}. He analyzed various form of structured sets in his \cite{schr74}, pervening to the definitions of {\it group}, {\it semi-group} and {\it loop}. More on this topic in my \cite{algebra} and \cite{gruppi}.} For this ground, reaching a satisfiying definition of {\it relation} is of fundamental importance, because on the concept of {\it relational structured theory} all Schr\"oder's mathematical work revolves. With this goal in sight, Schr\"oder devoted the third volume of his {\it Vorlesungen \"uber die Algebra der Logik} \cite{vorl3} to an analysis of the concept of {\it relation}.\\
\indent Because of some ambiguity on the side of Schr\"oder, these investigations culminate at the same time in a {\it language by signs} [Zeichensprache] or \emph{Pasigraphy}, able to express the main concepts of the exact sciences, and to free mathematics from the chains of the natural speech.\footnote{See \cite[leaf 2]{barton} or \cite[p.~75]{schr74}. Here Schr\"oder speaks of a language which {\it turns out to be more a language by signs than a language by words}.} Schr\"oder was not alone in envisaging an universal language in the 19th century. This search was typical of the period. After Latin ceased to be the {\it lingua franca} for scientists and humanists, the need of a substitute was urgent in order to facilitate the dialogue between people speaking different languages. This is not a trivial matter of human comunication (and understanding). I think to the language of music. Thanks to its universality is understandable by anyone. The same score can be performed by a German, by a Russian, by an Inuit (Eskimo), and this way the performers can interchange their interpretations.\footnote{It is my opinion that Joyce wrote his {\it Finnegangs Wake} in order that it could be approached the same way by people from different mother langues.}\\
\indent Unfortunately, things are a little more involved in natural languages, than in music. We cannot presuppose the knowledge of any language by a mathematician to make mathematics. Such requirement would get rid of those scholar without linguistic abilities. Today, the lingua franca of mathematics (and of many other disciplines as well) is English, but not in Schr\"oder's time. In the second half of the 19th century, many languages could pretend to be the new lingua franca; but preferring one language to another would caused that any nation supported own language with a sort of nationalism.\\
\indent But a discipline, as mathematics, which is famous for its universality could not accept the rule of {\it one} language over another. This was surely the ground backing up the search of a new lingua franca. I am obliged to stress again Schr\"oder's formalism. It would make no sense expressing formulas spoiled of any meaning whatever, by a natural language with all its references to an {\it interpreted} world. It is best forging a new language from scratch.\\
\indent Schr\"oder' supposed lingua franca was called by him {\it Pasigraphy}\footnote{From the ancient Greek: pan/pas- = to all + graphein = writing. According to wikipedia, a pasigraphy is a writing which must be understood by people from different linguistic areas. See http://de.wikipedia.org/wiki/Pasigrafie. Last visit the 2nd October 2013. The word 'pasigraphy' was coined, by Joseph de Maimieux (1753--1820) in his 1797 book, entitled {\it Pasigraphie ou premieres \'el\'ements du novel art-science d'ecrire et d'imprimer en une langue de mani\`ere \`a \^{e}tre lu et entendu dans toute autre langue sans traduction} \cite{vola3}. de Maimieux introduced a purely symbolic language, translatig any word belonging to natural language in a particular artificial sign. Schr\"oder made no reference at all for his pasigraphy, mentioning only the {\it Volap\"uk} \cite[p.~75]{barton}, an artificial language created by the German priest Johann Martin Schleyer (1831--1912) in 1879--1880 ca. Schleyer called his language {\it Volap\"uk}, which  means {\it international language} \cite[p.~369]{vola4}, being composed by vol- (= world, universe \cite[p.~368]{vola4}) and -p\"uk (= language, speech, tongue, dialect \cite[p.~273]{vola4}). I don't know from whom Schr\"oder borrowed the name 'pasigraphy'. What I may tell is that a certain Karl Obermair in 1864 founded a pasigraphical association (unfortunately I found no notices on him). Moses Pei\'c, the same year, distinguished between a pasi{\it graphy} (writing) and a pasi{\it logy} (spoken language) in his {\it System einer Universalsprache sowohl durch die Schrift (Pasigrafie), als auch durch die Laute (Pasilogie)}, published in Wien the 1864 \cite{vola2}. Pei\'c's pasigraphy is simply a vocabulary of numbers. To any word of natural language is assigned a number. This way, a discourse is a collection of different numbers. Personally, I consider de Maimieux' and Pei\'c's efforts a little too exotic. Perhaps, Schr\"oder knew these or similar work. I don't believe that he created the name of 'pasigraphy', from the ancient Greek, by himself. There are too many striking coincidences.} and was no other than the {\it Calculus of Relatives} considered as \emph{language}. Once laid down the fundamental laws governing this language (the so called {\it Fundamental Statements} [Festsetzungen])\footnote{For the most part, these statments are rules governing the construnction of well-formed formulas in the calculus of relatives. A couple of them express properties of well formed formulas in this calculus.}, and exposed the main problem of the calculus (i.e. the {\it Solution Problem}), Schr\"oder shows that his language is capable to express all mathematics. Really, he is conted in showing that Arithmetics can be expressed in the calculus of relations (or Pasigraphy)\footnote{I remember Schr\"oder's ambiguity on the role of the calculus of relatives: calculus and language at the same time.}, following Richard Dedekind, who in his masterwork {\it Was sind und was sollen die Zahlen?} stated:
\begin{quotation}
(\ldots) it appears as something self-evident and not new that every theorem of algebra and higher analysis, no matter how remote, can be expressed as a theorem about natural numbers -- a declaration I have heard repeatedly from the lips of Dirichlet.\footnote{\cite[p.~792]{dedekind}.}
\end{quotation}
Well, if any theorem whatever can be espressed in Arithmetic, it's only matter to found arithmetics to found mathematics. That explains Dedekind' and Schr\"oder's restriction to arithmetics. Then, Schr\"oder, translating Dedekind's {\it Theory of Chains}\footnote{The Theory of Chains was developed by Dedekind in order to found arithmetics.} in the calculus of relatives is showing that the calculus of relatives is not only able to express {\it some} fundamental mathematical concept, but {\it any} mathematical concept tout-court.\\
\indent The lecture of the third volume of the {\it Vorlesungen} devoted to this topic is the ninth. I tempted to say that this lecture has a virtuostic appeal. A sort of exercises to be performed, before to face the important mathematical questions of the time, more or less, as a piano player studies sets of exercises in order to perform a difficult score.\\
\indent One could questions my interpretation, quoting Schr\"oder himslef:
\begin{quotation}
The {\it ultimate} goal of the work [i.e. the translation of Dedekind's {\it Kettenlehre} in the calculus of relatives] is: to achieve a rigorous [streng] logical {\it definition} of the {\it relational} concept {\it number of-}, from which all sentences relating to this concept are to be derived in a pure deductive way.\footnote{\cite[pp.~349--350]{vorl3}. Translations from Schr\"oder's texts are mine, if not otherwise stated.}
\end{quotation}
I may reply to this important objection, that in this excerpt Schr\"oder is not declaring a foundational goal, but a {\it re-writing} one; i.e. how can the concept of {\it number} be translated in terms of binary relations? Obviously, it must be trasformed in a relation, that of {\it number of-}. This translation symbolizes in the calculus of relatives the concept of {\it number} from Dedekind.\footnote{For this "translation", see \cite[pp.~43--50]{bondoni}.}\\

\subsection{A little exemplification}\label{example}
I will show a couple of Dedekind's theorems translated in terms of relations. Let $a,b,c,$ be binary relations whatsoever. $1$ denote the universe of thought. The sign $;$ denotes the composition of relations. If we assume for a while that $f$ and $x$ are relations, $f;x$ is the relational translation of $f(x)$:\footnote{Obviously, in the following, $\mathfrak{D}$ stays for 'Dedekind'. Notice that Schr\"oder had a no clear idea of connectives. He wrote $(b\subseteq c) \subseteq (a;b \subseteq a;c)$ meaning $(b\subseteq c)\to (a;b \subseteq a;c)$.}
\begin{align}
&\mathfrak{D}\text{22}.\quad (b \subseteq c) \subseteq (a;b \subseteq a;c).\notag\\
&\mathfrak{D}\text{23}.\quad a;(b+c+\ldots) = a;b+a;c+\ldots \vert\vert\quad \mathfrak{D}\text{24}.\quad a;bc\ldots\subseteq a;b\cdot a;c\ldots\notag\\
&\mathfrak{D}\text{36}.\quad\text{Def.}\quad (a;b\subseteq b)= (\text{$a$ maps $b$ in itself})\notag\\
&\mathfrak{D}\text{37}.\quad\text{Def.}\quad  (a;b\subseteq b)= (\text{$b$ is a chain under $a$})\notag\\
&\mathfrak{D}\text{38}.\quad a;1\subseteq 1\notag\\
&\mathfrak{D}\text{39}.\quad (a;b\subseteq b) \subseteq (a;a;b\subseteq a;b)\notag
\end{align}

\noindent \cite[p.~354]{vorl3}. Roughly speaking, a chain is a relation which is closed under an application ($a;b \subseteq b$). Now, we state the correspondenting sentences by Dedekind. $A,B,C,S$ and $Z$ are systems (what today we call \emph{sets}), $A'$ is the image of $A$ under some function $\phi$, $\mathfrak{M}$ denotes the \emph{union of systems}, $\mathfrak{G}$ denotes the \emph{overlapping of systems}, and $\backepsilon$ the relation of inclusion:\\
\indent 22. Theorem. If $A \backepsilon B$, then $A'\backepsilon B'$.\\
\indent 23. Theorem. The image of $\mathfrak{M}$ ($A,B,C,\ldots$) is $\mathfrak{M}$($A',B',C',\ldots $).\\
\indent 24. Theorem. The image of every common part of $A,B,C,\ldots$, and therefore that of the intersection $\mathfrak{G}$ ($A,B,C,\ldots$) is part of $\mathfrak{G}$($A',B',C',\ldots$).\\
\indent 36. Definition. If $\phi$ is a similar\footnote{I.e. injective.} or dissimilar mapping of a system $S$, and $\phi (S)$ is a part of a system $Z$, then $\phi$ is said to be a mapping of $S$ into $Z$, and we say $S$ is mapped by $\phi$ into $Z$.\\
\indent 37. Definition. $K$ is called a \emph{chain} [Kette] when $K'\backepsilon K$.\\
\indent 38. Theorem. $S$ is a chain.\\
\indent 39. Theorem. The image $K'$ of a chain $K$ is a chain.\footnote{\cite[pp.~800--803]{dedekind}.}\\

I will not be polemic, but as it is evident from a comparison between Schr\"oder's list and Dedekind's one, we cannot speak of foundation at all. Does $\mathfrak{D}$22 found the theorem 22 by Dedekind? I don't believe. On the contrary, a comparison shows that Schr\"oder and Dedekind are speaking of the same matter with different \emph{languages}. I may only concede to my oppositors that a relation is more general than a function. But no more. The closure of a chain (theorem 37 by Dedekind) is perfectly mimicked by $\mathfrak{D}$37 in Schr\"oder. In any case, no one of the above theorems and definitions by Schr\"oder are more fundamental than the corresponding ones by Dedekind.\\
\indent But this is not all. Who can understimate the role of the principle of induction in Dedekind's work? The following is the principle in issue as laid down by Dedekind:\\

\indent 80. Theorem of complete induction (inference from $n$ to $n'$). In order to show that a theorem holds for all numbers $n$ of a chain $m_0$, it is sufficient to show,\\
\indent\quad $\rho$. that it holds for $n=m$, and\\
\indent\quad $\sigma$. that from the validity of the theorem for a number $n$ of the chain $m_0$ its validity for the following number $n'$ always follows.\footnote{\cite[p.~809]{dedekind}.}\\

\noindent This is the \emph{Gegensatz} in Schr\"oder:
\begin{equation}
\mathfrak{D}\text{59}.\quad\text{Theorem of complete induction}:\lbrace a;(a_0;b)c + b \subseteq c\rbrace \subseteq (a_0;b\subseteq c)\notag\footnote{\cite[p.~355]{vorl3}. Schr\"oder denotes a chain with a $0$ as subscript.}
\end{equation}
Analyzing the above sentence will carry us too far. It is sufficient to note that the request that $b\subseteq c$ corresponds to the first premise ($\rho$) in Dedekind, and that $a;(a_0;b)c \subseteq c$ corresponds to the second premise above ($\sigma$). Ergo, the conclusion ($a_0;b\subseteq c$) holds too. To use Schr\"oder's own words:
\begin{quotation}
We can rephrase in natural language the sentence $\mathfrak{D}$59 -- considered for the time being uniquely as a theorem on binary relations --  in the following way: \emph{To prove, that the $a$-chain\footnote{Speaking of $a$-chain, Schr\"oder made explicit the map under which a relation is closed; in this case, $a$. I remind that I am speaking of maps, but these maps are binary relatives, not functions.} of a relative $b$ is included in a third relation $c$, we need only} to make two things; i.e. it is sufficient \emph{to show}:\\
\indent first, \emph{that $b$ is included in $c$},\\
\indent second, \emph{that also $a$-image of every [ordered] couple of elements belonging to the $a$-chain of $b$, which are included in $c$, is also included in $c$.}\\
In other words,  $a_0;b$ must be part of $c$ [$a_0;b\subseteq c$], as soon as $b$ is part of $c$ [$b\subseteq c$] and the $a$-image of any common elements to  $a_0;b$ and $c$ [$a;(a_0;b)c$] is also part of $c$.\footnote{\cite[p.~367]{vorl3}. If $a;(a_0;b)c\subseteq c$ and $b\subseteq c$, then $a_0;b \subseteq c$, which is $\mathfrak{D}$59.}
\end{quotation}

We may rephrase the principle of induction also in the following manner: $(b\subseteq c\wedge (a;b\subseteq c\to a;(a;b)\subseteq c))\to (a_0;b\subseteq c).$\footnote{In Appendix A, I will prove this principle.}
As it is easily to see, also in this case we have only a translation of a principle stated in terms of sets in an analogue (albeit more general) in the calculus of relatives.\\
\indent Only one question reamins open. We asserted that the concept of relation generalizes that of function or set, because these are only a particular cases of relations. One could be tempted to say that it is that generality to constitue a fundament; i.e. the theorems in the calculus of relations found their pendants in set theory because they are more general, more uncompassing. It could be, but nowhere Schr\"oder states such position, limiting to shed light on the more perspicuity and elegance of his symbolic calculus. In other words, for Schr\"oder is a matter of rethoric, not of founding. In fact, let give voice again to Schr\"oder:
\begin{quotation}
With this [i.e. the symbolic language of the relatives] the reader has at his disposal a \emph{key} to \textbf{translate} one representation of the chain theory in the other. As one can see, our method of symbolize [Bezeichnungsweise] is the \emph{most expressive}. (\ldots) our representation of the chain theory is so no way inferior in \emph{clarity} to any other -- neither to that of a such master of precision and concision, who is its author [i.e. Dedekind].\footnote{\cite[p.~353]{vorl3}. The emboldening is mine.}
\end{quotation}
And in a paper of 1895:
\begin{quotation}
I will here not insist, that in our discipline [i.e. the calculus of relatives] we succeeded in condensing even more the sentences of such a master of concision [i.e. Dedekind] (\ldots).\footnote{\cite[p.~157]{schr95}.}
\end{quotation}
Such excerpts show the rethorical possibilities of the calcuus of relatives. It is in the context that it arises the goal to translate Dedekind's concept of number in a more elegant language. Schr\"oder want exhibit the power of his theory:
\begin{quotation}
Meanwhile, hoping not too late, I will now face the task to \emph{include} Dedekind's ``Theory of Chain'' in the edifice [Lehrgeb\"aude] of our discipline in order to give a proof of its [i.e. of the calculus of relatives] power. The result (\ldots) will be suitable to demonstrate for the first time the value of our discipline.\footnote{\cite[p.~346]{vorl3}.}
\end{quotation}
Schr\"oder could be not more clear: translating Dedekind's theory of chain is an usefull exercise to show the symbolic value of the calculus of relatives. Value which is more rethorical than computational, as the quotations above on the expressivity of the calculus of relations are to witness. Of course, one may object that the tentative to condense long formulas will not always have as a result more understandable formulas. This is another matter. Schr\"oder believed in his language.

\subsubsection{A turn in Schr\"oder's thought}
Given the importance of this question I will insist upon. Schr\"oder's goal is to re-write the concept of {\it number} in his calculus of relatives now regarded as a pasigraphy, i.e. an universal and symbolic language. Given the ambiguity of the calculus in issue (theory and language), it is evident that this task waves between foundationalism and language; but Schr\"oder had no real foundational interest, mantaining that foundationalism had already a definitive solution, that proposed by Dedekind. No more efforts were necessary: Dedekind gave the definitive answer.\footnote{As already noted, Schr\"oder did not realize that Dedekind introduced the definition of {\it set} of natural numbers, and not the definition of a {\it single} number $n$.}\\
\indent If until this point of the \emph{Vorlesungen} Schr\"oder was engaged in investigating the concept of 'relation', because from this concept depends its structural philosophy of mathematics, now there is a shift in this thought: \emph{can the calculus of relatives not only solve some theoretical problem, but also serve as a symbolic language?} For Schr\"oder the answer is adfirmative, and indeed in the following sections of the \emph{Vorlsegunen} Schr\"oder focuses on translating some set-theoretical pivotal notions in the calculus of relatives: set, function, etc.

\subsection{Brady on Schr\"oder's foundationalism}
My interpretation relies on considering Schr\"oder a mathematician and not a logician as usually. As a matter of fact, Schr\"oder knew of logic by authors who were eminent mathematicians, as Boole and de Morgan, both engaged in analysis (derivation and integration) and not from philosophers.\footnote{With the utmost probability, Schr\"oder knew the work of Robert Grassmann, by the brother of the late, Hermann Grassmann, a \emph{mathematician} engaged in Vector Calculus.} Furthermore, many scholars considers as epiphanic the short autobiographical sketch in \cite{barton}, which I am not sure was written by Schr\"oder, at least in its entirity. Really, that sketch is a sort of publicity for which Schr\"oder payed. In any case, despite my own interpretation, I repute correct to give voice to some scholars mantaining that Schr\"oder wrote his third volume of {\it Vorlesungen}\footnote{In particular, the ninth lecture.} in order {\it to found} mathematics.  I will exemplify this interpretation relying on Brady's book \cite{brady}. Obviously this point of view is common (in some case, only partly) to the eminent historian Volker Peckhaus, to Risto Villko and to Javier Legris. For the respective positions of these three scholars, I refer to the bibliography. In this place, for sake of clarity I will take in consideration only the work of Geraldine Brady, leaving aside further declinations.
\begin{quotation}
[The {\it Vorlesungen \"uber die Algebra der Logik} offer] the first exposition of abstract lattice theory,\footnote{Despite a radicate tradition, it was not Schr\"oder to prove that non any lattice is distributive, being the proof by J.~L\"uroth. I don't understand why Schr\"oder don't make the name of his friend, originating a misunderstanding.} the first exposition of Dedekind's theory of chains after Dedekind, the most comprehensive development of the calculus of relations, and \textbf{a treatment of the foundations of mathematics in relation calculus} that L\"owenheim in 1940 still thougth was as reasonable as set theory.\footnote{\cite[p.~143]{brady}. The emboldening is mine.}
\end{quotation}
And some page below,
\begin{quotation}
Schr\"oder translates Dedekind's set-theoretic treatment of chains line-by-line into the second-intentional calculus of relatives. With this, Schr\"oder shows that the second-intentional theory of relatives is sufficient to develop number theory.\footnote{\cite[p.~158]{brady}.}
\end{quotation}
The antecedent of this quotation is right: as a matter of fact, as seen before, Schr\"oder \emph{translated} Dedekind's theory of chains in his calculus of relatives, but Geraldine Brady draws from this antecedent a {\it false} consequent. She does not take in account that Schr\"oder's calculus of relations was \emph{both} a calculus and a symbolic language. When we speak of \emph{translating}, the calculus of relatives as language is meant. Schr\"oder limited himself to re-write the Kettenlehre in his calculus of relatives. Is it sufficient such re-writing to speak of {\it foundationalism}? If I translate the word {\it death} in German as {\it Tod}, I am not founding the English concept of death in the German one. It is only a question of translating a string of symbols in another string of symbols. No foundation is required.\\
\indent But Brady insists, referring to Leopold L\"owenheim and Alfred Tarski:
\begin{quotation}
The Peirce-Schr\"oder theme that higher intentional relative calculus can be a full foundation for mathematics recurs twice in later mathematical history. First, L\"owenheim \cite{low} made the claim that the relative calculus was just as suitable for a foundation of mathematics as set theory. Second, the theme of {\it Set Theory without Variables} \cite{tarski2} of Tarski and Givant (1987) is that a form of binary relation calculus is adequate as a foundation for all of mathematics, and uses no variables.\footnote{\cite[p.~159]{brady}. For this theme, see \cite[Chapter 4, pp.~67--92]{bondoni}.}
\end{quotation}
That both L\"owenheim and Tarski had foundational goals is manifest; that Schr\"oder was a source of inspiration for them is also true. What is false, is that Schr\"oder too was a foundationalist. He inspired foundationalists, without being himself a foundationalist.\\
In the same page Brady quotes C.S.~Peirce:
\begin{quotation}
The nearest approach to a logical analysis of mathematical reasoning that has ever been made was Schr\"oder's statements (\ldots) in a logical algebra of my invention, of Dedekind's reasoning (\ldots) concerning the foundations of arithmetics.\footnote{\cite[p.~344]{peirce} quoted in \cite[p.~159]{brady}.}
\end{quotation}
But it is not all. Brady, referring to the Lecture in which Schr\"oder faces the Kettentheorie, states:
\begin{quotation}
It seems likely that the purpose of this lecture was to show that the most delicate piece of foundations work thus far in the history of mathematics could be carried out neatly in the calculus of relatives.\footnote{\cite[p.~296]{brady}.}
\end{quotation}
What Brady fails to appreciate is that the calculus of relatives for Schr\"oder was not only a theory but also a language. Schr\"oder was engaged in finding a symbolic and universal language for mathematics.

\subsection{Not to found: that is the question}
Then, if the ninth lecture has not a foundational r\^{o}le, that of pervening to a definition of {\it set of natural numbers}, what is its real meaning? We must keep in mind that for Schr\"oder the calculus of relations was not only a calculus, but also a language in which expressing the main concepts of the exact sciences, as Legris rightly states:
\begin{quotation}
Algebra of relatives is considered both as a {\it universal language}, and as a theory on which any scientific science can be {\it founded}.\footnote{\cite[p.~243]{legris}.}
\end{quotation}
I cannot but agree with Legris on the duplice nature of the calculus of relatives, language and calculus (or theory) at the same time; notwithstanding, I don't believe that the calculus of relatives had for Schr\"oder a foundational character, as this paper is aiming to prove. For the German mathematician, the calculus of relatives was only a lattice of formulas devoid of meaning, and so capable of many interpretations (models). And from this fact, Schr\"oder' {\it structuralism} arises.\\
\indent It can be sound strange that the calculus of relations is {\it both} a language and a calculus, but we must keep in mind that for Schr\"oder a language is devoid of any interpretation. It is not Frege's {\it Conceptual Notation} which has a canonical interpretation \cite{revlegris}; i.e. it refers to only one interpretation. From this point of view, Frege's language is nearer to the natural languages than Schr\"oder's one.\\
\indent No canonical interpretation is presupposed by Schr\"oder, the words of his language being only \emph{formal} well-formed formulas. His language tells nothing. It needs a model to become informative. It is not by chance that Schr\"oder called his language a {\it language by signs}.\\
\indent Schr\"oder's calculus of relatives is completely formal, being a lattice of well-formed formula which are strings of inkspots (signs) on paper. It is this formality which the calculus of relatives and an universal language share. Calculus and language are two sides of the same coin.

\section{On pasigraphy}\label{found2}
We can now pass to the second task of this paper: summing up some evidence to showing that Schr\"oder's pasigraphy was connaturate to his author and not an hobby which Schr\"oder cultivates in his last years. First of all, I repute interesting to question {\it why Schr\"oder abandoned the old name of {\it Language by Sign}, in favour of {\it Pasigraphy}}? Why did he need a name for his calculus? Notice that Schr\"oder was not the unique to envisage an artificial and universal language in the 19th century; I think to the Esperanto, to the Volap\"uk, to Frege {\it conceptual notation}, to Peano {\it latino sine flexione}, etc.\\
\indent The Karlsruhe mathematician, probably, choose the name of {\it pasigraphy}, because of its neutrality. {\it Pasigraphy}, as noted above, means {\it universal language}, and Schr\"oder with his calculus of relatives aimed just to it, to a universal language for mathematics.\\ 
\indent For Schr\"oder the pasigraphy\footnote{We must say "Schr\"oder's pasigraphy", because his pasigraphy was not the unique. Volap\"uk, Esperanto or other similar linguistic efforts were all a pasigraphy, an universal language.} was the best possible language for mathematics and in order to show this, he translated Dedekind's Chain Theory in his language. The ninth lecture of the third volume of the {\it Vorlesungen} is so not devoted to foundational aims, but to exhibit the power and elegance of the calculus of relations, as mantained above:
\begin{quotation}
I will here not insist, that in our discipline [i.e. in the calculus of relatives], we succeeded in condensing even more the sentences of such a master of concision [i.e. Dedekind](\ldots).\footnote{\cite[p.~157]{schr95}.}
\end{quotation}
From this quotation is manifest the care Schr\"oder gave to in shaping his language. Any symbol of it is carved with the utmost attention, as it is the case for the symbols for sum and product of relations:
\begin{quotation}
Because the non commutative behaviour of the sum between relatives, I shaped [gestaltet] the plus-sign not symmetrically; Peirce, instead, managed with the erected Cross as in death notices. For similar reasons, I choose for the relative multiplication the {\it semicolon}, because this is a not symmetrical sign, adapt to represent a not symmetrical composition (\ldots)\footnote{\cite[p.~33]{vorl3}. See the table in the Appendix B.}
\end{quotation}
In other words, from the shape of the sign it must be evident his r\^{o}le. In fact, the non symmetrical sign of "$;$" denotes and suggests a {\it not} commutative ({\it not} symmetrical) operation. We could say that Schr\"oder's pasigraphy is really an {\it Ideography}, a {\it pictorial} symbolism, a language by {\it pictures}.\footnote{For contemporary similar efforts, see at least \cite{coocke}.} The meaning of a sign is suggested by its picture. Take for example the symbol $\dag$. Its picture suggests the death, being the image of the cross on which our Lord Jesu Christ died.\\
\indent The first time Schr\"oder uses the word "Pasigraphie" is in his 1890 delivered prolusion \cite{schr90}, six years before his famous contribute at the {\it First International Mathematical Congress} held in Z\"urich:\footnote{See \cite{schr98}. Schr\"oder was appointed {\it rector} at the Karlsruhe University for the academic year 1890--1891. He held his oath the 3rd November 1890, see \cite[leaf 28]{diener}. Notice that Schr\"oder was rector for only one year. In this sense, the traditional opposition between a Peirce without any academical engagement and Schr\"oder who succeed in making career in the academical milieu must be revisited. The appointment of 1890, probably was due to a lack of a better scholar. Significantly, the lapse of time in which Schr\"oder lead the Kalrsruhe University is very small.}
\begin{quotation}
Such system of signs [Bezeichnung] is, once extended to the entire field of the objects of thought, in opposition to the signs of the [natural] words, which are more or less equivalent from the point of view of the content they represent, a typical language of concepts, a conceptual notation [Begriffsschrift]\footnote{Please, notice the expression!}, and in contrast to the various languages used by normal people, a general language of the thing, a Pasigraphy or a Universal Language.\footnote{\cite[p.~16]{schr90}. We know from before that {\it Pasigraphy} is synonimous of {\it universal language}.}
\end{quotation}
Notice what Schr\"oder is stating: while {\it the signs of the words of a natural language are more or less equivalent from the point of view of the content they represent}, the signs of the calculus of relatives manifest their content by their picture. It is not so in many natural languages, where the words have a denotation which is independent from their graphical appearance.\\
\indent Take again the symbol $\dag$. We said that it reminds the reader to the death, being a pictorial symbol. Take now the word 'death'; nothing in its graphical form suggests that 'death' denotes the death. It is only a matter of convention. Obviously, there is no unique symbol to indicate ideographically the death, but the question is not on the univocity of a denotation, but in the ability to denote by a picture.\\
%
%
\indent Any way, compare the last quotation from the {\it Vorlesungen} with the incipit of \cite{vola3}:
\begin{quotation}
The word PASIGRAPHY is composed by two Greek words, {\it to all}, and GRAPHO, {\it I am writing}. To write even to whom who does not understand any language, \textbf{by a writing which is a picture of the thought}, represented by different syllables, that is what we call PASIGRAPHY.\footnote{\cite[p.~1]{vola3}. The emboldening and the translations from this book are mine.}
\end{quotation}
The idea underlying de Maimieux's efforts is to create a language which is a picture of the things, in plain agreement with Schr\"oder who just spoke of a {\it general language of things} and in his {\it Vorlesungen} shaped his signs according to the content they must convey. For this reason, a semicolon suggests visually the concept of non symmetry. In other words, we grasp that a such operation is not a symmetrical/commutative one by the real picture of its sign on paper. It is not requested any further knowledge. But let us quote againg from \cite{vola3}:
\begin{quotation}
(\ldots) a text, hand-written or printed [using the pasigraphy], can be read and understood in many languages, as the arithmetics in ciphers; then the characters of chemistry and of music are the same way intelligible from Petersburg to Malta, from Madrid to Pera, from London and Paris to Philadelphia or to the Bourbons.\footnote{\cite[ivi]{vola3}.}
\end{quotation}
What is interesting in these two last quotations by de Maimieux is the definiton of pasigraphy as language by signs. So we found the reason why Schr\"oder introduced in his vocabulary the word 'pasigraphy': it denotes an universal language by pictures (signs). There is not a phasigraphical phase in Schr\"oder's thought as maintained from some scholars: from 1873 onwards, Schr\"oder spoke of a language by signs.\footnote{See \cite[pp.~350--352]{LL1}. In the 1873 {\it Lehrbuch der Arithmetik und Algebra}, Schr\"oder speaks of signs as concrete objects attached to papers, comparing them to {\it mushrooms} (sic!) \cite[pp.~16--17]{schr73}.} Now, he attached a more {\it popular} name to it.\\  
\indent Obviously, I don't believe that Schr\"oder knew the work of de Maimieux or that of Pei\'c. He cited only the Volap\"uk as source of inspiration; but it is not to be excluded that Schr\"oder found in Schleyer's books some reference to previous linguistic efforts named {\it pasigraphy}.\\
\indent Finally, we must not forget Schr\"oder's interest in learning diverse languages:
\begin{quotation}
At the age of eight, due largely to his grandfather's encouragement, he could read Latin. He later acquired proficiency, to varying degrees, in French, English, Italian, Spanish, and Russian. This linguistic ability eventually allowed him to correspond with Poretskii and other Russian logicians in Russian, with Peirce, Ladd-Franklin, John Venn and sometimes even his fellow German Paul Carus in English. His ability in Italian allowed him to read the work of, and later correspond with, Peano and Padoa. These linguistic abilities thus placed him in an important position in logic of the late 19th century, which was increasingly becoming an \textbf{international discipline}, with major works in English, German, and Italian (or Peano's {\it latino sine flexione}).\footnote{\cite[p.~120]{dipert2}. The emboldening is mine.}
\end{quotation}
Unfortunately, this quotation seems exhibiting a Schr\"oder whose only merit was to organize the calculus of relatives and to mantain links between logicians from diverse linguistic areas. It is not so. I quoted this long excerpt by Dipert only because it stress Schr\"oder's love for speaking diverse language. I don't share Dipert's overall interpretation of Schr\"oder.\\
\indent What it is interesting is that Schr\"oder's pasigraphy has many roots: on one side, it is the result of Schr\"oder genuine love for {\it language}, on the other side, it is the consequence of his engagement in abstract fields of mathematics, as algebra; finally, it was a necessity to overspread worldwide his work. Dipert stressed the importance of language in Schr\"oder's everyday life; Volker Peckhaus, Risto Villko and Legris stated the connection between abstract algebra and a pasigraphy. I inserted Schr\"oder in the broader context of the search for an artificial universal language.

\subsection{Set-Theory}
All this discourse aimed to confute a pretese foundationalism in Schr\"oder. The third Volume of the {\it Vorlesungen} is mainly devoted to non-foundational problems; the most part of the third volume of {\it Lectures on the Algebra of Logic} handles the algebraic solution problem (5th Lecture) and from the 10th Lecture onwards, set-theory (definitions of set, of finite and infinite set; definition of function, injective and bijective, equivalence, and so on). On these topics the last Schr\"oder will ponder.\\
\indent Often, one reads that Schr\"oder's work on relation was interrupted by his death. Schr\"oder finished the first part of the 3rd volume of the {\it Lectures} in 1895. He will die only in 1902. In these seven years he puts aside the theory of relations, considered in itself, to study set-theoretic problems expressed in his new forged language of relatives. To be precise, the following are the last papers by Schr\"oder:
\begin{enumerate}
\item Note \"uber die Algebra der bin\"aren Relative (1895) \cite{schr95},\footnote{This paper is a very short summary of the ninth lecture from the third volume of the {\it Vorlesungen}.}
\item \"Uber Pasigraphie, ihren gegenw\"artigen Stand und die pasigraphische Bewegung in Italien (1897) \cite{schr98},
\item Ueber zwei Definitionen der Endlichkeit und G.~Cantor'sche S\"atze (1898) \cite{schrcan2},
\item Die selbst\"andige Definition der M\"achtigkeiten $0,1,2,3$ und die explizite Gleichzahligkeitsbedingung (1898) \cite{schrcan1},
\item On Pasigraphy. Its Present State and the Pasigraphic Movement in Italy (1899) \cite{schr99},
\item \"Uber G.~Cantorsche S\"atze (1901) \cite{schr1901a},\footnote{In this short note, Schr\"oder states the conditions under which two set can be said {\it equivalent}.}
\item Sur une extension de l'id\'ee d'ordre (1901) \cite{schr01},
\item Ernst Schr\"oder (short autobiography) (1901) \cite{barton}.
\end{enumerate}
We must not neclegt \cite{lovett} who describes Schr\"oder's talk at the International Congress of Philosophy, held in Paris in 1900. According to Lovett, that occasion Schr\"oder spoke of relations and well-ordening. Then, if we erase from the above listing the first and the last item, we may appreciate in what the last Schr\"oder was engaged: set theory and pasigraphy. No hint to a foundationalism is present.
\newpage

\section*{Appendix A}
In this short appendix, we will prove the principle of induction in the calculus of relatives. In order to be as clear as possible, I will use the principle in issue in my own formulation, reminding that it is equivalent to Schr\"oder's one. Then, we must demonstrate that:
\begin{equation}
(b\subseteq c\wedge (a;b\subseteq c\to a;a;b\subseteq c))\to (a_0;b\subseteq c).
\end{equation}
By hypothesis, the generator of the chain $a_0;b$ belongs to $c$:
\begin{equation}\label{uno1}
b\subseteq c
\end{equation}
Then, by theorem $\mathfrak{D}$36,\footnote{See above in the section \ref{example}.}, by \eqref{uno1}, and by transitivity:
\begin{equation}\label{due}
a;b\subseteq c.
\end{equation}
From the same theorem $\mathfrak{D}$36 we can easily infer that $a;a;b\subseteq a;b$; from this, \eqref{due}, and transitivity:
\begin{equation}\label{tre}
a;a;b\subseteq c.
\end{equation}
Iterating $n$ times the process, using $\mathfrak{D}$36, our hypothesis $b\subseteq c$ and transitivity, we obtain:
\begin{equation}\label{enne}
\overbrace{a;a;\ldots;a;}^{\text{$n$ times}}b\subseteq c.
\end{equation}
Let us simplify \eqref{enne} using Schr\"oder's own definition,
\begin{equation}\notag
\text{Def.}\quad\overbrace{a;a;\ldots;a;}^{\text{$n$ times}}b \equiv a_{00};b.\footnote{\cite[p.~340]{vorl3}.}
\end{equation}
The definition above enable us to compactify \eqref{enne}:
\begin{equation}\label{enne2}
a_{00};b\subseteq c.
\end{equation}
From \eqref{uno1} and \eqref{enne2}, we can state:
\begin{equation}\label{enne3}
b\wedge a_{00};b\subseteq c.
\end{equation}
A theorem by Dedekind, $\mathfrak{D}$58, states:
\begin{equation}\label{dedekind11}
a_0;b \equiv b \wedge a_{00};b.\footnote{\cite[p.~375]{vorl3}.}
\end{equation}
Finally, from \eqref{enne3} and \eqref{dedekind11} it follows:
\begin{equation}\label{fine}
a_0;b\subseteq c. \qquad\square
\end{equation}

\subsection{Commentary}
This proof is highly interesting for many reasons. Let us see why. First of all, notice that in theorem $\mathfrak{D}$36 we may omit the reference to the domain $b$ which is closed under the map in issue. The result is surprising:
\begin{equation}\notag
a;a\subseteq a.\footnote{\cite[p.~337]{vorl3}.}
\end{equation}
This is the condition which a relative $a$ must satisfy in order to be transitive. If we take in account of this fact, it will be manifest that in the calculus of relatives we don't need the concept of chain. It is sufficient to require that the relative $a$ be transitive. In other words, this proof gets rid of the concept of set, of similar mapping and of chain. The major part of Dedekind's work is useless in the calculus of relations. As a matter of fact, to prove the principle of induction, which leads us to the concept of set of natural numbers, in the calculus of relatives is sufficient that the relation under which a relative is closed be transitive, and furthermore that we can express a chain by its generators plus the iteration of the mapping in question. That $a_0;b$ is equivalent to $\overbrace{a;a;\ldots;a;}^{\text{$n$ times}}b$ is shown already at page 326 of \cite{vorl3}, just before Schr\"oder's investigations on transitivity.\\
\indent Taking in account that both the transitivity and the shortcut for $\overbrace{a;a;\ldots;a;}^{\text{$n$ times}}b$ are in the \textbf{eight} lecture and not in the ninth (devoted to the chain theory), our proof is saying that Schr\"oder did not need to translate Dedekind's chain theory in his one to found mathematics. It could accomplish this task, before introducing the chain theory. Ergo, we obtained a further rationale to deny the foundational goals in Schr\"oder. If he would found mathematics, he could do it already in the eight lecture, but he did not. As a matter of fact, Schr\"oder aimed not to found mathematics, but to cast light on the power of his calculus, which is another question.\\
\indent Finally, all the stuff we employed in order to prove the principle of induction is part of a work on the \emph{Solution Problem}. Is all this not sufficient to persuade my opponents form whom I am \emph{not a serious scholar}, that Schr\"oder had a mathematical point of view in analyzing the calculus of relations?\footnote{Obviously, in this appendix I used the calculus of relations as calculus and not as a symbolic language.} 
\newpage

\section*{Appendix B}
In the following a list of the symbols employed by Schr\"oder is provided. For typographical reason I used the symbols $\subseteq$ and $\cdot'$ instead of the original ones. Inside square brackets a modern pendant is introduced.
\begin{table}[h]
\begin{center}\renewcommand{\arraystretch}{1.3}
\begin{tabular}{|c|l|}
\hline
\large{\textbf{Schr\"oder' Symbols}} & \large{\textbf{Meaning}}\\ \hline
$a,b,c,\ldots$ & binary relations [$R,S,T,\ldots$]\\ \hline
$i,j,k,\ldots$ & as subscripts, individual variables [$x,y,z,\ldots$]\\ \hline
$0$ & empty relation\\ \hline
$1$ & universe of thought, or universal relation\\ \hline
$1'$ & diagonal [$\delta_{ij}$; $\delta$-Kronecker operator]\\ \hline
$0'$ &\emph{anti}-diagonal\\ \hline
$=$ & equal [or equivalence]\\ \hline
$\subseteq$ & improper inclusion [also implication]\\ \hline
$+$ & sum [union, but also $\vee$ and sometimes $\wedge$]\\ \hline
$\cdot$ or simple juxtaposition & times [intersection, but also $\wedge$]\\ \hline
$;$ & composition [$\circ$]\\ \hline
$\cdot'$ & percean sum [$\bullet$]\\ \hline
$\prod$ & universal first-order quantifier\\ \hline
$\sum$ & existential first-order quantifier\\ \hline
\end{tabular}
\end{center}
\end{table}

\nocite{bondoni}
\nocite{revlegris}
\nocite{peck1}
\nocite{peck2}
\nocite{peck3}
\nocite{peck4}
\nocite{vola1}
\nocite{vola2}
\nocite{vola3}
\nocite{schr99}
\nocite{coocke2}
\nocite{brady}
\nocite{vilkko}
\nocite{alex2}

\bibliographystyle{amsalpha}
\bibliography{hamlet}

\end{document}